PRESENTATIONS OF FINITELY GENERATED GROUPS}

\documentstyle[twoside,epic,eepic,12pt]{article}
\title{\bf THE NORMALIZED CYCLOMATIC QUOTIENT ASSOCIATED WITH PRESENTATIONS OF
FINITELY GENERATED GROUPS}

\author
{Amnon Rosenmann \\
\vspace {-2 mm}
\small
Dept. of Math. \& Computer Science \\
\vspace {-2 mm}
\small
Ben-Gurion University \\
\vspace {-2 mm}
\small
Beer-Sheva, Israel \\
\vspace {-2 mm}
\small
aro@black.bgu.ac.il
}

\date{}

\begin{document}
\maketitle

\newtheorem{claim}{Claim}[section]
\newtheorem{remark}[claim]{Remark}
\newtheorem{example}[claim]{Example}
\newtheorem{definition}[claim]{Definition}
\newtheorem{theorem}[claim]{Theorem}
\newtheorem{lemma}[claim]{Lemma}
\newtheorem{proposition}[claim]{Proposition}
\newtheorem{corollary}[claim]{Corollary}
\newtheorem{conjecture}[claim]{Conjecture}

\newcommand{\GAia} {\mbox{$G_1^{\alpha_1}$}}
\newcommand{\GAib} {\mbox{$G_2^{\alpha_2}$}}
\newcommand{\GAir} {\mbox{$G_r^{\alpha_r}$}}
\newcommand{\GAi} {\mbox{$G_i^{\alpha_i}$}}
\newcommand{\GAj} {\mbox{$G_j^{\alpha_j}$}}
\newcommand{\GAa} {\mbox{$G^{\alpha_1}$}}
\newcommand{\GAb} {\mbox{$G^{\alpha_2}$}}
\newcommand{\GA} {\mbox{$G^{\alpha}$}}
\newcommand{\HAia} {\mbox{$H_1^{\alpha_1}$}}
\newcommand{\HAi} {\mbox{$H_i^{\alpha_i}$}}
\newcommand{\HA} {\mbox{$H^{\alpha}$}}

\newcommand{\stFFM} {\mbox{$\scriptstyle {\mathcal F}^*$}}
\newcommand{\stFF} {\mbox{$\scriptstyle {\mathcal F}$}}
\newcommand{\FFM} {\mbox{${\mathcal F}^*$}}
\newcommand{\FF} {\mbox{$\mathcal F$}}
\newcommand{\stCFM} {\mbox{$\scriptstyle {\mathcal C} {\mathcal F}^*$}}
\newcommand{\stCF} {\mbox{$\scriptstyle {\mathcal C} {\mathcal F}$}}
\newcommand{\CFM} {\mbox{${\mathcal C} {\mathcal F}^*$}}
\newcommand{\CF} {\mbox{${\mathcal C} {\mathcal F}$}}
\newcommand{\CC} {\mbox{$\mathcal C$}}
\newcommand{\DD} {\mbox{$\mathcal D$}}
\newcommand{\CS} {\mbox{$\mathcal S$}}
\newcommand{\TT} {\mbox{$\mathcal T$}}

\newcommand{\dout} {\mbox{$\scriptstyle out$}}
\newcommand{\ds} {\mbox{$\scriptscriptstyle S$}}
\newcommand{\dx} {\mbox{$\scriptscriptstyle X$}}
\newcommand{\dt} {\mbox{$\scriptscriptstyle \mathcal{T}$}}
\newcommand{\dd} {\mbox{$\scriptscriptstyle D$}}
\newcommand{\bdt} {\mbox{$\scriptscriptstyle \overline{\mathcal T}$}}
\newcommand{\dfhB} {\mbox{$\scriptscriptstyle F/H_2$}}
\newcommand{\dfhA} {\mbox{$\scriptscriptstyle F/H_1$}}
\newcommand{\dhB} {\mbox{$\scriptscriptstyle H_2$}}
\newcommand{\dhA} {\mbox{$\scriptscriptstyle H_1$}}

\newcommand{\dfh} {\mbox{$\scriptscriptstyle F/H$}}
\newcommand{\dg} {\mbox{$\scriptscriptstyle G$}}
\newcommand{\dh} {\mbox{$\scriptscriptstyle H$}}
\newcommand{\df} {\mbox{$\scriptscriptstyle F$}}
\newcommand{\dfa} {\mbox{$\scriptscriptstyle F/A$}}
\newcommand{\da} {\mbox{$\scriptscriptstyle A$}}
\newcommand{\dc} {\mbox{$\scriptscriptstyle C$}}

\newcommand{\lrw}{\mbox{$\longrightarrow$}}
\newcommand{\Lrw}{\mbox{$\Longrightarrow$}}
\newcommand{\Llrw}{\mbox{$\Longleftrightarrow$}}
\newcommand{\SemiDP}{\mbox{$> \! \! \! \lhd$}}

%\begin{abstract}
%\input{abs}
%\end{abstract}

\baselineskip=12pt  %line can be omitted if changed to 12pt##
\section{Introduction}
Given the Cayley graph of a finitely generated group $G$, with respect to a
presentation $\GA$ with $n$ generators,
the quotient of the rank of the fundamental group of subgraphs of the
Cayley graph by the cardinality of the set of
vertices of the subgraphs gives rise to the definition of the normalized
cyclomatic quotient $\hat \Xi(\GA)$.
The asymptotic behavior of this quotient is similar to the
asymptotic behavior of the quotient of the cardinality of the boundary
of the subgraph by the cardinality of the subgraph. Using F{\o}lner's
criterion for amenability one gets that $\hat \Xi(\GA)$ vanishes for infinite
groups if and only if they are amenable. When $G$ is finite then $\hat \Xi(\GA)
= 1/|G|$, where $|G|$ is the cardinality of $G$,
and when $G$ is non-amenable then $1-n \leq \hat \Xi(\GA)<0$,
with $\hat \Xi(\GA) = 1-n$ if and only if $G$ is free of rank $n$.
Thus we see that on special cases $\hat \Xi(\GA)$ takes the values of the
Euler characteristic of $G$. Most of the paper is concerned with
formulae for the value of $\hat \Xi(\GA)$ with respect to that
of subgroups and factor groups, and with respect to the decomposition of
the group into direct product and free product. It is not surprising that
some of the formulae and bounds we get for $\hat \Xi(\GA)$ are similar to
those given for the spectral radius of symmetric random walks on the
graph of $\GA$ (see \cite{Kes} and also \cite{Gri0} for the connection of
the spectral radius with the ``growth-exponent''),
but this is not always the case (see e.g. the remark after
Theorem~\ref{thC40}). The difference is that we do not explore the graph
associated with a given presentation in a ``breadth-first'' manner but rather in
a ``depth-first'' manner, proceeding in the ``best'' possible way.
In the last section of the paper we define and touch very briefly
the {\em balanced} cyclomatic
quotient, which is defined on
concentric balls in the graph. This definition is
related to the growth of $G$.

Throughout the paper we assume that when we are given presentations
$\GAi$ of groups $G_i$ and $H$ is a group defined through the $G_i$,
then $H$ gets the natural {\em induced} presentation $\HA$. Thus, if $H$
is a factor of $G$ then $\HA$ has the same generating set as that of the
presentation $\GA$ of $G$ with all the relations of $\GA$ holding
also in $\HA$. Or if $H=G_1 * G_2$ then the generators and relators of
$\HA$ are the union of those of the $\GAi$, assuming that the generating
sets of the $\GAi$ are disjoint. Similarly for direct products (with
the appropriate commutators as extra relators), etc. We do not try,
however, to be too precise with regard to the use we make of the notation $\GA$
or $\GAia$.

The formulae we give for direct and free products give upper and lower
bounds for $\hat \Xi(\GA)$ in the following sense. Suppose that
$\GAi = < X_i \mid R_i > $, $i=1,2$, and $\GA = < X \mid R > $, with
$X = X_1 \cup X_2$, $R \supseteq R_1 \cup R_2$ and both unions disjoint,
so that the natural maps $\phi_{i} : G_{i} \rightarrow G$, $i=1,2$,
are injective. Then
\begin{equation}
\hat{\Xi}(\GAia * \GAib) \leq \hat \Xi(\GA) \leq
\hat{\Xi}(\GAia \times \GAib).
\end{equation}
The left inequality follows from Theorem~\ref{thC40} since the Cayley
graph associated with $\GA$ is a quotient of that associated with
$\GAia * \GAib$. The right inequality appears in the proof of
Theorem~\ref{thC20}. These bounds hold for such structures as
semi-direct products, amalgamated products, HNN extensions, etc.
One can also try and get exact formulae for the structures mentioned above,
at least in special cases, but we do not get much into it in the present
paper. \\

We use the following terminology and notation on graphs.
The set of vertices of a graph $\mathcal G$ is denoted by $V({\mathcal G})$
and the set of vertices by $E({\mathcal G})$. A {\em path}
in $\mathcal G$ is a sequence $v_0,e_1,v_1,e_2, \ldots$,
$v_i \in V({\mathcal G})$, $e_i \in E({\mathcal G})$, such that $e_i$
starts at the vertex $v_{i-1}$ and terminates at $v_i$. The
length of a path $v_0,e_1,v_1,e_2, \ldots, v_n$ is $n$.
A {\em simple path} is a path in which the vertices along it are distinct,
except possibly for the first and last one, in which case it is a
{\em simple closed path} or a {\em simple circuit}.
We assume that each path is {\em reduced}, i.e. it is not homotopic
to a shorter one when the initial and terminal vertices are kept fixed.

If ${\mathcal H} \subseteq {\mathcal G}$, i.e. $\mathcal H$ is a collection
of vertices
and edges of the graph $\mathcal G$, then
we denote by $< {\mathcal H} >$ the subgraph {\em generated} by $\mathcal H$.
It is the smallest
subgraph of $\mathcal G$ which contains $\mathcal H$. That is, we add
to $\mathcal H$ the
endpoint vertices of all the edges in $\mathcal H$. On the other hand, the
subgraph of $\mathcal G$ {\em induced} by $\mathcal H$ is the one
whose vertices are
those of $\mathcal H$ and whose edges are all the edges which join
these vertices
in $\mathcal G$. An induced subgraph is a subgraph which is induced by some
${\mathcal H} \subseteq {\mathcal G}$. If ${\mathcal H}_1, {\mathcal H}_2
\subseteq {\mathcal G}$
then ${\mathcal H}_1 - {\mathcal H}_2$ is the collection of
vertices $V({\mathcal H}_1) - V({\mathcal H}_2)$ and edges
$E({\mathcal H}_1) - E({\mathcal H}_2)$,
and it does not necessarily form a subgraph of $\mathcal G$, even when
${\mathcal H}_1$
and ${\mathcal H}_2$ are subgraphs of $\mathcal G$. The {\em boundary}
of the subgraph $\mathcal H$ of $\mathcal G$ is
$\partial {\mathcal H} =  {\mathcal H} \cap <{\mathcal G} - {\mathcal H}>$,
and its {\em interior} is $\dot{\mathcal  H} = {\mathcal H} - \partial
{\mathcal H}$.
The {\em outer boundary} of $\mathcal H$ (in $\mathcal G$) is the set
of vertices
of ${\mathcal G} - {\mathcal H}$ which are adjacent to $\mathcal H$ in
$\mathcal G$.
Assume now that each edge of $\mathcal G$ is labeled with some $x \in X$
in one direction and with $x^{-1} \in X^{-1}$ in the other direction.
Then we define $E^X_{out}({\mathcal H})$ to be the set of edges of
${\mathcal G} -
{\mathcal H}$ whose initial vertices with respect to the directions $X$ are
in $\mathcal H$.

Finally, let $\beta_0({\mathcal H}) = |V({\mathcal H})|$,
let $\beta_1({\mathcal H}) = |E({\mathcal H})|$, and let $\alpha({\mathcal H})
= |\pi_0({\mathcal H})|$ be the number of the (connected)
components of $\mathcal H$.
\section{The Normalized Cyclomatic Quotient}
Let $\GA =<X \mid R>$ be a presentation of a group $G$, with
$X=\{x_1, \ldots, x_n \}$, and let $\mathcal G$ be the associated Cayley graph.
If $H$ is the normal closure of $R$ in $F=<X>$ then it is shown in
\cite{Ros4} that the growth function of $G$ is equivalent
to the {\em rank-growth} $rk_{\dh}$ of $H$. The rank-growth is defined by
\begin{equation}
rk_{\dh} (i) = \mbox{rank} (H_i),
\end{equation}
where $H_i$ is the subgroup of $F$ generated by the elements of length
$\leq i$ (with respect to $X \cup X^{-1}$). We notice that $H_i$
is the fundamental group of the subgraph of $\mathcal G$ of all paths
starting at $1$ of length $\leq i$.
Thus, there exists an {\em exhausting chain} $({\mathcal G}'_i)$
in $\mathcal G$,
i.e. a sequence ${\mathcal G}'_1 \subseteq {\mathcal G}'_2 \subseteq
{\mathcal G}'_3
\subseteq \cdots$ of subgraphs of $\mathcal G$ whose union is $\mathcal G$,
with all the ${\mathcal G}'_i$ connected and finite, such that
the growth of the function $\gamma_1(i) = \mbox{rank}(\pi_1({\mathcal G}'_i))$
is equivalent to the growth of the function $\gamma_2(i)=
|V({\mathcal G}'_i)|$.
We are now interested in the asymptotic behavior of the quotient
$\gamma_1(i) / \gamma_2(i)$, which is related, as we will see, to the
quotient $| V(\partial {\mathcal G}'_i) | / |V({\mathcal G}'_i)|$. The latter
quotient and its analogs are known and widely studied objects in diverse
areas of Mathematics (see e.g. the survey \cite{Li}). By F{\o}lner's
criterion (see
\cite{Gro}) the group $G$ is amenable if and only if there exists an
exhausting chain $({\mathcal G}'_i)$ of finite connected subgraphs of
$\mathcal G$ such that
\begin{equation}
\limsup_{i \rightarrow \infty} \frac{\beta_0(\partial {\mathcal G}_i)}
{\beta_0({\mathcal G}_i)} = 0
\label{eqC64}
\end{equation}
(recall that a group $G$ is amenable if there exists
an invariant mean on $B(G)$, the space of all bounded complex-valued
functions on $G$ with the sup norm $\parallel f \parallel_{\infty}$, see
\cite{Gre}). We remark that in F{\o}lner's criterion one can consider
as well disconnected subgraphs or boundaries of any fixed width
$k$. Notice that (\ref{eqC64}) implies that if $G$ is non-amenable then
it has exponential growth.

Let us denote by $\beta_2({\mathcal G}')$ the cyclomatic number of
${\mathcal G}'$,
i.e. the sum of the values of $\mbox{rank} (\pi_1({\mathcal H}'))$ over all
the components ${\mathcal H}'$ of the subgraph ${\mathcal G}'$
(the notation $\beta_2({\mathcal G}')$ refers to the number of $2$-cells of
an associated $2$-dimensional complex). Let us also use the following
notation:
\begin{eqnarray}
\xi({\mathcal G}') &=& \frac{\beta_2({\mathcal G}')}{\beta_0({\mathcal G}')},
\\ \mu({\mathcal G}') &=& \frac{\beta_1({\mathcal G}')+1}
{\beta_0({\mathcal G}')}.
\end{eqnarray}
We denote by $\FF({\mathcal G}), \FFM({\mathcal G}), \CF({\mathcal G}),
\CFM({\mathcal G})$
respectively the sets of finite, non-trivial finite, connected finite,
non-trivial connected finite subgraphs of $\mathcal G$.
Here a graph is non-trivial if it contains more than one vertex.
\begin{definition}
{\em
If $C=({\mathcal G}_1 \subseteq {\mathcal G}_2 \subseteq \cdots)$,
${\mathcal G}_i \in \CF({\mathcal G})$, is an exhausting chain, let
\begin{equation}
\eta_{\dc} = \limsup_{i \rightarrow \infty} \xi({\mathcal G}_i).
\end{equation}
Then we define the {\em cyclomatic quotient} of $\GA$ by
\begin{equation}
\Xi(\GA) = \sup_{\dc} \eta_{\dc}, \ \ \ \ \ \
\mbox{$C$ an exhausting chain},
\end{equation}
and the {\em normalized cyclomatic quotient} of $\GA$ by
\begin{equation}
\hat \Xi(\GA) =1-n+ \Xi(\GA).
\label{eqC68}
\end{equation}

The following are equivalent definitions of $ \hat \Xi(\GA)$.
\begin{eqnarray}
\hat \Xi(\GA) &=& 1-n+ \sup_{{\mathcal G}' \in \stFF({\mathcal G})}
\xi({\mathcal G}'),
\label{eqC70} \\
\hat \Xi(\GA) &=& \sup_{{\mathcal G}' \in \stCF({\mathcal G})}
\frac{1-|E^X_{out}({\mathcal G}')|} {\beta_0({\mathcal G}')},
\label{eqC80} \\
\hat \Xi(\GA) &=& -n + \sup_{{\mathcal G}' \in \stCF({\mathcal G})}
\mu({\mathcal G}'),
\label{eqC82} \\
\hat \Xi(\GA) &=& \sup_S
\left( \frac{1}{|S|}+\sum_{j=1}^{n} \left( \frac{|Sx_j \cap S|}
{|S|} - 1 \right) \right),
\ \ \ \  \mbox{$S \subseteq G$ finite}.
\label{eqC83}
\end{eqnarray}
}
\end{definition}
In case $G$ is infinite then we have
\begin{equation}
\begin{array}{ll}
\hat \Xi(\GA) = - \inf_{{\mathcal G}' \in \stCF({\mathcal G})}
\frac{|E^X_{out}({\mathcal G}')|} {\beta_0({\mathcal G}')}, &
\label{eqC80a} \\
\hat \Xi(\GA) = -n + \sup_{{\mathcal G}' \in \stCF({\mathcal G})}
\frac{\beta_1({\mathcal G}')}{\beta_0({\mathcal G}')}, &  \\
%\label{eqC82a}
\hat \Xi(\GA) = \sup_S
\sum_{j=1}^{n} \left( \frac{|Sx_j \cap S|}
{|S|} - 1 \right),
& \mbox{$S \subseteq G$ finite}.
%\label{eqC83a}
\end{array}
\end{equation}
Clearly definition (\ref{eqC70}) gives at least the same value as
definition (\ref{eqC68}). To see that these definitions are equivalent
we need to show that for every ${\mathcal G}', {\mathcal G}'' \in
\FF({\mathcal G})$
there exits ${\mathcal H} \in \CF({\mathcal G})$ such that ${\mathcal G}''
\subseteq {\mathcal H}$
and $\xi({\mathcal G}') \leq \xi({\mathcal H})$. But this follows by
Lemma~\ref{lmC20}
and by the fact that given any subgraphs ${\mathcal G}_1, {\mathcal G}_2 \in
\FF({\mathcal G})$, we can cover ${\mathcal G}_2$ by the induced subgraph of
translates of ${\mathcal G}_1$.
\begin{lemma}
Let ${\mathcal G}' \in \FF({\mathcal G})$ satisfy $\xi({\mathcal G}') \geq
\xi({\mathcal H}')$
for every ${\mathcal H}' \subseteq {\mathcal G}'$. Let ${\mathcal G}'' \in
\FF({\mathcal G})$
such that $\xi({\mathcal G}'') \geq \xi({\mathcal G}')$, and let ${\mathcal H} =
{\mathcal G}' \cup {\mathcal G}''$. Then $\xi({\mathcal H}) \geq
\xi({\mathcal G}')$.
\label{lmC20}
\end{lemma}
{\em Proof}.
Clearly
\begin{equation}
\beta_0({\mathcal H}) = \beta_0({\mathcal G}') + \beta_0({\mathcal G}'') -
\beta_0({\mathcal G}' \cap {\mathcal G}''),
\end{equation}
while
\begin{equation}
\beta_2({\mathcal H}) \geq \beta_2({\mathcal G}') + \beta_2({\mathcal G}'') -
\beta_2({\mathcal G}' \cap {\mathcal G}'').
\end{equation}
Then by simple calculation we get that
\begin{equation}
\xi({\mathcal H}) \geq \xi({\mathcal G}')
\label{eqC69}
\end{equation}
after clearing denominators in the quotients $\xi({\mathcal H})$ and
$\xi({\mathcal G}')$.
\hfill $\Box$ \\

Definition (\ref{eqC80}) follows by
\begin{equation}
\beta_2({\mathcal G}') =
1 + (n - 1)\beta_0({\mathcal G}') - |E^X_{out}({\mathcal G}')|, \ \ \ \ \
{\mathcal G}' \in \CF({\mathcal G}).
\label{eqC74}
\end{equation}

Definition (\ref{eqC82}) follows by
\begin{equation}
\beta_2({\mathcal G}') =
\alpha({\mathcal G}') - \beta_0({\mathcal G}') + \beta_1({\mathcal G}')
\end{equation}
($\alpha({\mathcal G}')$ is the number of components of ${\mathcal G}'$).

Finally, for definition (\ref{eqC83}) we notice that if
${\mathcal G}'$ is an induced subgraph of ${\mathcal G}$ and
that $S \subseteq G$ is the set of elements of $G$ which correspond to
the vertices of ${\mathcal G}'$ then
\begin{equation}
\beta_1({\mathcal G}') = \sum_{j=1}^{n} |Sx_j \cap S|.
\label{eqC81}
\end{equation}

The slightly simplified expressions (\ref{eqC80a}) when $G$ is infinite
are due to the fact that ${\mathcal G}'$ may be chosen as large as we wish.
\begin{proposition}
Let $\GA$ be a presentation of a group $G$ with $n$ generators and let
${\mathcal G}$ be the associated Cayley graph. Then $1-n \leq \hat
\Xi(\GA) \leq 1$.
Moreover,
\begin{description}
\item[(i)] if $G$ is amenable then $\hat \Xi(\GA)=1 / |G|$, where $1 / |G|$
is defined to be $0$ if $|G|=\infty$;
\item[(ii)] if $G$ is non-amenable then $1-n \leq \hat \Xi(\GA) < 0$,
with $\hat \Xi(\GA)= 1-n$ if and only if $G$ is free of rank $n \geq 2$.
\end{description}
\label{prC77}
\end{proposition}
{\em Proof}.
If $G$ is finite then any exhausting chain of ${\mathcal G}$ stabilizes on
${\mathcal G}$. Then by (\ref{eqC80})
\begin{equation}
\hat \Xi(\GA) = \frac{1-|E^X_{out}({\mathcal G})|}{\beta_0({\mathcal G})}=
\frac{1}{|G|}.
\end{equation}
In fact, we see that for every proper subgraph ${\mathcal G}'$ of ${\mathcal G}$
\begin{equation}
\frac{1-|E^X_{out}({\mathcal G}')|}{\beta_0({\mathcal G}')} \leq 0.
\end{equation}

Since for every ${\mathcal G}' \in \FF({\mathcal G})$,
$|E^X_{out}({\mathcal G}')|$ is of the same order as
$\beta_0(\partial {\mathcal G}')$, then by using F{\o}lner's criterion we obtain
from (\ref{eqC80}) that $\hat \Xi(\GA)=0$ when $G$ is infinite amenable.

When $G$ is non-amenable then there exists $c>0$ such
that for every ${\mathcal G}' \in \FF({\mathcal G})$,
$|E^X_{out}({\mathcal G}')| / \beta_0({\mathcal G}') > c$.
Then by (\ref{eqC80}) $\hat \Xi(\GA) \leq -c$, by letting $\beta_0({\mathcal G}'
)
\rightarrow \infty$. On the other hand, $\hat \Xi(\GA)> 1-n$ when there exits at
least one circuit in ${\mathcal G}$. When ${\mathcal G}$ contains no circuits
then $G$ is free of rank $n$ and $\hat \Xi(\GA)= 1-n$.
\hfill $\Box$ \\

If we use (\ref{eqC82}) for the definition of $\hat \Xi(\GA)$ then by
Proposition~\ref{prC77} we get the following criterion for amenability:
$G$ is amenable if and only if for every $\epsilon > 0$ there exists
${\mathcal G}' \in \FF({\mathcal G})$ such that
\begin{equation}
\frac{\beta_1({\mathcal G}')}{\beta_0({\mathcal G}')} > n-\epsilon.
\end{equation}
In other words, $G$ is amenable if and only if
for every $\epsilon > 0$ there is a finite subset $S$ of $G$ such that
$|Sx_j \cap S| / |S| > 1-\epsilon$ for each
generator $x_j$. This is easily seen to be equivalent to F{\o}lner's
criterion for amenability (\cite{Fol}):
for every $\epsilon > 0$ and every $w_1,\ldots,w_r \in G$ there is a
finite subset $S$ of $G$ such that
$|Sw_i \cap S| / |S| > 1-\epsilon$ for each $i$. (This also shows that
a subgroup of an amenable group is amenable.)
\\

The value of $\hat \Xi(\GA)$ depends on the presentation $\GA$ of $G$.
Let us look at the following example. Let $\GAa$ be
the presentation of the free group $G$ of rank $2$ with generators $x,y$.
Then $\hat \Xi(\GAa)=1-2=-1$. Let $\GAb$ be obtained from $\GAa$ by the
Tietze transformation of adding a new generator $x'$ and a relation $x'=w$.
If $w = x^k$ for some integer $k$ then $\hat \Xi(\GAb)=
\hat \Xi(\GAa)$, as can be seen from the chain of subgraphs whose vertices
consist of increasing powers of $x$. On the other hand, if $w=xy$ then
the only simple circuits we get are triangles of the form
$x'y^{-1}x^{-1}$
(or cyclic permutations of it) and by forming an increasing chain of
subgraphs the best we can do is adding each time $2$ new vertices
and obtaining a new circuit. Therefore we get that $\hat \Xi(\GAb)=(1-3)+
1/2 = -3/2 < \hat \Xi(\GAa)$.
\begin{proposition}
Let $\GA=< X \mid R>$, $X=\{x_1, \ldots, x_n\}$.
\begin{description}
\item[(i)] If $\GAa=<X \cup \{x'\} \mid R, \ x'=w(X) >$, where $w(X) \in <X>$,
then $\hat \Xi(\GA)-1 \leq \hat \Xi(\GAa) \leq \hat \Xi(\GA)$, with
$\hat \Xi(\GAa) = \hat \Xi(\GA)$
if $w=1$.
\item[(ii)] If $\GAb=<X \cup \{x'_1, \ldots, x'_n\} \mid R, \ x'_i=x_i, \
i=1, \ldots, n>$ then $\hat \Xi(\GAb)=2 \hat \Xi(\GA) - 1/|G|$.
\end{description}
\label{prC20}
\end{proposition}
{\em Proof}.
{\bf (i)} The Cayley graph of $G$ with respect to $\GAa$ is obtained from the
Cayley graph of $G$ with respect to $\GA$ by adding the edges in the
direction $x'$ from each vertex $v$ to the vertex $vw$. Thus we can increase
$\beta_1({\mathcal G}')$ by at most $\beta_0({\mathcal G}')$. The result then
follows from (\ref{eqC82}).
When $G$ is (finite or infinite) amenable
then $\hat \Xi(\GAa)=\hat \Xi(\GA)$ by Proposition~\ref{prC77}. One can
also see it
directly from (\ref{eqC80}) by considering the ``thickening'' of ${\mathcal G}'$
to ${\mathcal G}''$ by adding to it the outer $d$-boundary, where $d=l(w)$,
and noticing that $E^{X'}_{out}({\mathcal G}'') / \beta_0({\mathcal G}'')
\rightarrow 0$,
where $X'=X \cup \{x'\}$. When $w=1$ we can make sure that the edges going-out
in directions $X$ are the same as those going-out in directions $X'$ and thus ob
tain
$\hat \Xi(\GAa)=\hat \Xi(\GA)$.

{\bf (ii)} When $G$ is finite then $\hat \Xi(\GAb)=\hat \Xi(\GA)=1/|G|$ by
Proposition~\ref{prC77}. When $G$ is infinite then since the number of
out-going edges is doubled we get that $\hat \Xi(\GAb)=2 \hat \Xi(\GA)$ as
$\hat \Xi(\GAb) = - \inf_{{\mathcal G}'} |E^X_{out}({\mathcal G}')| /
\beta_0({\mathcal G}')$, ${\mathcal G}' \in \CF({\mathcal G})$.
\hfill $\Box$ \\

We notice that by Proposition~\ref{prC77} and Proposition~\ref{prC20}~(ii)
we get that $\hat \Xi(\GA)$ is independent of the presentation if and
only if $G$ is amenable. \\

When $G$ is finite then $\xi$ has, of course, a maximum, and we have
seen that the maximum is achieved only at the whole graph ${\mathcal G}$ of the
presentation. We showed also that in general (for finite and infinite
groups), for every proper subgraph ${\mathcal H}$ of ${\mathcal G}$ there is a
subgraph ${\mathcal H}'$ which properly contains ${\mathcal H}$ and such that
$\xi({\mathcal H}') \geq \xi({\mathcal H})$. We then ask more: does $\xi$ have a
maximum in case the group is infinite? It is quite clear that when $G$
is the free product of a finite group and a free group, with a
``natural'' presentation so that we have generators of the free part
which are not involved in any relator, $\xi$ does have a maximum
because no circuit involves the generators of the free factor (for
more details on the value of $\xi$ on free products see
section~\ref{secFP}). Next we will show that in fact this is the only
situation where $\xi$ has a maximum.

Given a presentation $\GA=<X \mid R>$, with $X=\{x_1, \ldots, x_n \}$
and $R$ not empty, let $c(x_i)$ be the length of a shortest
relator in which $x_i$ appears (i.e. the shortest (reduced) circuit in
the Cayley graph which contains the edge $x_i$) or $0$ if $x_i$ does not
appear in any relator. Then let $c(\GA)=\max_i(c(x_i))$.
\begin{theorem}
Let $\GA=<X \mid R>$, $|X|=n$. Then one of the following holds.
\begin{description}
\item[Case 1:] $\hat \Xi(\GA)=1-n+\xi({\mathcal G}')$ for some ${\mathcal G}' \i
n
\FF({\mathcal G})$.
Let $X' \subseteq X$ be the set of
labels of the edges not going-out of ${\mathcal G}'$. Then
\begin{description}
\item[(i)] if $X'=X$ then ${\mathcal G}'={\mathcal G}$ and $G$ is finite;
\item[(ii)] if $X'=\emptyset$ then $G$ is the free group on $X$;
\item[(iii)] if $|X'|=k$, $1 \leq k \leq n-1$, then $V({\mathcal G}')$ is the
union of left cosets of the finite subgroup $H=Gp(X')$ (which may be trivial)
and $G=H*F$ where $F$ is the free group on $X-X'$.
\end{description}
\item[Case 2:] $\xi$ does not have a maximum on $\FF({\mathcal G})$.
Then for every ${\mathcal G}' \in \FF({\mathcal G})$
\begin{equation}
\hat \Xi(\GA) \geq 1-n+\xi({\mathcal G}')+\frac{1}{(c(\GA)-1)
\beta_0({\mathcal G}')}.
\label{eqC05}
\end{equation}
\end{description}
\label{thC10}
\end{theorem}
{\em Proof}.
Case 1: it suffices to show that whenever there is a relator involving
an edge going-out of ${\mathcal G}'$ then $\hat \Xi(\GA) > \xi({\mathcal G}')$.
For suppose to the contrary that $\hat \Xi(\GA) = \xi({\mathcal G}')$ and
that $\lambda=v_0,e_1,v_1, e_2, \ldots,e_m,v_0$ is a simple circuit
starting at $v_0 \in \partial{\mathcal G}'$ and that $e_1 \in
E^X_{out}({\mathcal G}')$.
Note that $m \geq 2$
since $\xi({\mathcal G}')$ is maximal. Then by Lemma~\ref{lmC20}
we can add translates of ${\mathcal G}''$ along the vertices
of $\lambda$ which are not in ${\mathcal G}'$ so that the resulting subgraph
${\mathcal G}''$ will satisfy $\xi({\mathcal G}'') \geq \xi({\mathcal G}')$.
Moreover, if
${\mathcal G}''$ does not contain the edge $e_1$ then by adding this edge to
${\mathcal G}''$ we
increase $\xi$ - in contradiction to assumption. To see that this is
indeed possible, let $y_k \in X \cup X^{-1}$ be
the label of $e_k$. Then we need to show that for every $k$ for which
$v_k \notin V({\mathcal G}')$,
there exists $u_k \in V({\mathcal G}')$ such that $u_k w_k^{-1} \notin
V({\mathcal G}')$, where $w_k= y_1 \cdots y_k$ (here we look at the vertices
as the group elements  they represent). But this follows by the assumption
that $V({\mathcal G}')$ is not mapped to itself by the map $g \mapsto g w_k$.

Case 2: let ${\mathcal H}'_0 \in \FF({\mathcal G})$ and assume that
$\xi({\mathcal H}'_0)
\geq \xi({\mathcal H}')$ for every ${\mathcal H}' \subseteq {\mathcal H}'_0$,
or otherwise we will take such a subgraph ${\mathcal H}' \subseteq {\mathcal H}'
_0$
and get a better result.
Then, as we have shown above, there exists a subgraph ${\mathcal H}'_1$,
which is a union of translates of ${\mathcal H}'_0$,
with $\beta_0({\mathcal H}'_1) \leq  c(\GA) \beta_0({\mathcal H}'_0)$ and with
\begin{equation}
\xi({\mathcal H}'_1) \geq \xi({\mathcal H}'_0) + \frac{1}{c(\GA)
\beta_0({\mathcal H}'_0)}.
\end{equation}
The same process can now be carried out with ${\mathcal H}'_1$ and so on,
obtaining a sequence ${\mathcal H}'_i$ of subgraphs satisfying
\begin{equation}
\xi({\mathcal H}'_i) \geq \xi({\mathcal H}'_{i-1}) + \frac{1}{c(\GA)^i
\beta_0({\mathcal H}'_0)}.
\end{equation}
Thus, passing to the limit, we get that
\begin{equation}
\hat \Xi(\GA) \geq 1-n+\xi({\mathcal H}'_0)+\frac{1}{(c(\GA)-1)
\beta_0({\mathcal H}'_0)}.
\end{equation}
\hfill $\Box$ \\

We remark that in case 2 of the above theorem we get in particular, by
taking ${\mathcal G}'$ to be the trivial subgraph, that
\begin{equation}
\hat \Xi(\GA) \geq 1-n+\frac{1}{(c(\GA)-1)}
\end{equation}
and that there exists an exhausting chain $({\mathcal H}'_i)$, ${\mathcal H}'_i
\in \CF({\mathcal G})$, such that
\begin{equation}
\limsup_{i \rightarrow \infty}\frac{|E^X_{out}({\mathcal H}'_i)|}
{\beta_0({\mathcal H}'_i)} \leq  n-1-\frac{1}{(c(\GA)-1)}.
\end{equation}
\section{Subgroups and Factor Groups}
When $G_2$ is a homomorphic image of $G$, with the presentation
$\GAib$ induced by the presentation $\GA$, then the Cayley graph associated
with $\GAib$ may be regarded as a quotient of the Cayley graph associated
with $\GA$, which implies, as expected, that
$\hat \Xi(\GAib) \geq \hat \Xi(\GA)$. In fact we have the following.
\begin{theorem}
Let $G_1$ be a normal subgroup of $G$ and let
$G_2=G/G_1$ with the presentation $\GAib$ induced by $\GA$. Then
\begin{equation}
\hat \Xi(\GA) \leq \hat \Xi(\GAib) - \left( \frac{1}{|G_2|} - \frac{1}{|G|}
\right),
\label{eqC118}
\end{equation}
with equality holding if $G_1$ is amenable.
\label{thC40}
\end{theorem}
{\em Proof}.
We may exclude the cases where $G_1$ or $G_2$ are trivial.
When $G$ is finite the result follows by Proposition~\ref{prC77} since
$\hat \Xi(\GA) = 1 / |G|$.
So assume $G$ is infinite. We will prove first the inequality in
(\ref{eqC118}). Again by
Proposition~\ref{prC77} the result is clear when $G_2$ is finite.
So we further assume that $G_2$ is infinite.
Let $p : {\mathcal G} \rightarrow {\mathcal G}_2$ be the covering map from
the Cayley graph of $\GA$ onto that of $\GAib$.
If ${\mathcal G}' \in \CF({\mathcal G})$ then one can decompose ${\mathcal G}'$,
regarded as
a topological space, into a finite number of
subspaces $\bar{\mathcal E}'_i$, say $i=1, \ldots,r$, where each
$\bar{\mathcal E}'_i$ is in a different sheet, that is the map
\begin{equation}
p|_{\bar{\mathcal E}'_i} : \bar{\mathcal E}'_i \rightarrow {\mathcal E}'_i
\end{equation}
is bijective and continuous but not necessarily a homeomorphism.
Each ${\mathcal E}'_i$ is a subgraph of ${\mathcal G}'_2 = p({\mathcal G}')$.
Then $\beta_0({\mathcal G}') =
\sum_{i=0}^r \beta_0({\mathcal E}'_i)$, but $\beta_2({\mathcal G}')
\leq \sum_{i=0}^r \beta_2({\mathcal E}'_i)$, since having a simple circuit
of ${\mathcal G}'$ which lies
in $k$ sheets results in increasing $\sum_{i=0}^r \beta_2({\mathcal E}'_i)$
by $k$. Therefore $\xi({\mathcal G}') \leq \xi({\mathcal G}'_2)$. Since this
is true for any ${\mathcal G}' \in \CF({\mathcal G})$, we have
$\hat \Xi(\GA) \leq \hat \Xi(\GAib)$.

Suppose now that $G_1$ is amenable.
Then we need to show that $\hat \Xi(\GA) =
\hat \Xi(\GAib)$ in case $G_2$ is infinite, or that $G$ is amenable in
case $G_2$ is finite.
Let $X$ be the generating set of $\GA$ and let
${\mathcal G}, {\mathcal G}_2$ be the Cayley graphs associated with $\GA,\GAib$
respectively. Given $\epsilon > 0$, let ${\mathcal G}'_2 \in \CF({\mathcal G}_2)
$
such that
\begin{equation}
\frac{1-|E^X_{out}({\mathcal G}'_2)|} {\beta_0({\mathcal G}'_2)} > \hat \Xi(\GAi
b) -
\frac{\epsilon}{3}.
\end{equation}
Assume also that ${\mathcal G}'_2$ is induced, contains the vertex $1$, and
that ${\mathcal G}'_2={\mathcal G}_2$ in case $G_2$ is finite and otherwise
$\beta_0({\mathcal G}'_2) > 3 / \epsilon$. Let ${\mathcal T}'_2 \in
\CF({\mathcal G}_2)$
be a spanning tree of ${\mathcal G}'_2$. Each vertex of ${\mathcal T}'_2$ is the
n
assigned a specific element of $G$ (of which we make use in (\ref{eqC128})
below). ${\mathcal T}'_2$ is embedded as a tree ${\mathcal T}' \in
\CF({\mathcal G})$,
${\mathcal T}' \subseteq p^{-1}({\mathcal T}'_2)$ ($p$ the covering map) with
the same vertex and edge labels. Then we take ${\mathcal G}' \in
\CF({\mathcal G})$
to be the subgraph induced by ${\mathcal T}'$. Let $H_1$ be the subgroup
of $G_1$ generated by the set $Y$ consisting of the non-trivial elements
\begin{equation}
y_{v,x} = vx(p(vx))^{-1} \neq 1, \ \ \ \ \ \ \ v \in V({\mathcal T}'), \ x \in X
,
\ p(vx) \in V(T').
\label{eqC128}
\end{equation}
If $Y$ is empty we take $H_1$ to be the trivial group.
Let ${\mathcal H}_1$ be the Cayley graph of $H_1$ with respect to $Y$,
%If $G_1$ is finite let ${\mathcal H}'_1 = {\mathcal G}_1$, a Cayley graph
%of $G_1$,
and let
${\mathcal H}'_1 \in \CF({\mathcal H}_1)$ with
%$\beta_0({\mathcal H}'_1) > 4 / \epsilon$
\begin{equation}
\frac{1-|E^Y_{out}({\mathcal H}'_1)|} {\beta_0({\mathcal H}'_1)} >
-\frac{\epsilon \beta_0({\mathcal G}'_2)}{3}.
\end{equation}
Such a subgraph exists by the amenability of $H_1$.
Let ${\mathcal G}''
\in \FF({\mathcal G})$ be the subgraph induced by (the disjoint union)
$\bigcup_{g \in V(\scriptstyle {\mathcal H}'_1)} g {\mathcal G}'$
(we look here at $g \in V({\mathcal H}'_1)$ as an element of $G$ by
writing each $y \in Y$ with the generators $X$ of $G$).
Let now $e \in E^X_{out}({\mathcal G}'')$ be an
edge labeled with $x \in X$ and starting at $gv \in V({\mathcal G}'')$,
$g \in V({\mathcal H}'_1)$, $v \in V({\mathcal T}')$. Then either $p(e) \in
E^X_{out}({\mathcal G}'_2)$, or else $p(e)$ joins $v=p(gv)$ and $u=p(gu)$,
for some $u \in V({\mathcal T}')$. Then
\begin{equation}
gvx = gy_{v,x} p(vx) = gy_{v,x}u.
\end{equation}
That is, $e$ is the unique edge corresponding to an edge $e' \in
E^Y_{out}({\mathcal H}'_1)$ that starts at $g$ and
is in direction $y_{v,x}$, and this correspondence is 1-1. Then we have
\begin{eqnarray}
\frac{1 - |E^X_{out}({\mathcal G}'')|} {\beta_0({\mathcal G}'')} &=&
\frac{1-(\beta_0({\mathcal H}'_1)|E^X_{out}({\mathcal G}'_2)| +
|E^Y_{out}({\mathcal H}'_1)|)} {\beta_0({\mathcal H}'_1)\beta_0({\mathcal G}'_2)
} \\
&=&
\frac{1-|E^X_{out}({\mathcal G}'_2)|} {\beta_0({\mathcal G}'_2)} +
\frac{1-|E^Y_{out}({\mathcal H}'_1)|} {\beta_0({\mathcal H}'_1)
\beta_0({\mathcal G}'_2)} -
\frac{1}{\beta_0({\mathcal G}'_2)} \nonumber \\
&>&  \hat \Xi(\GAib) - \left(\frac{1}{|G_2|} - \frac{1}{|G|} \right)
- \epsilon, \nonumber
\end{eqnarray}
since $G$ is infinite and $|E^X_{out}({\mathcal G}'_2)| = 0$ if $G_2$
is finite.
That is, $\hat \Xi(\GA) \geq \hat \Xi(\GAib)$ if $G_2$ is infinite, and
$\hat \Xi(\GA) \geq 0$ if $G_2$ is finite. By the the inequalities in the
other directions - these are equalities.
\hfill $\Box$
\begin{corollary}
If $\GA=\GAia \SemiDP \GAib$ and $G_1$ is amenable then
\begin{equation}
\hat \Xi(\GA) = \hat \Xi(\GAib) - \left( \frac{1}{|G_2|} - \frac{1}{|G|}
\right).
\end{equation}
\label{crC30}
\end{corollary}
{\em Proof}.
By Proposition~\ref{prC20} we may assume that the presentation $\GAib$
is induced by the presentation $\GA$, by adding the generators $g_i$
of $\GAia$ and the relations $g_{i}=1$ for every $i$. Then the result
follows immediately by Theorem~\ref{thC40}.
\hfill $\Box$ \\

Another corollary to Theorem~{\ref{thC40} is the known fact that if both
$G_1$ and $G_2$ are amenable then $G$ is also amenable.

We remark that $G_1$ may be non-amenable but still $\hat \Xi(\GA) =
\hat \Xi(\GAib)$. For example, let $G=H*H*K$ where $H$ is a $2$-generated
finite group and $K$ is free of rank $2$. Let $K_1$ be a normal subgroup
of $K$ such that $K / K_1 \simeq H$, let $G_1$ be the normal closure
of $K_1$ in $G$, and let $G_2 = G / G_1$.
Then $G_2 \simeq H*H*H$ and by Corollary~\ref{crC40}~(i),
$\hat \Xi(G) = \hat \Xi(G_2)$ although
$G_1$ is non-amenable. This situation does not happen when considering
the spectral radii $R, R_2$ associated with symmetric random walks on
$G,G_2$ respectively, where $R = R_2$ if and only if $G_1$ is amenable
(see \cite{Kes}, Theorem~2).
\begin{definition}
{\em
Let $\GA$ be a presentation of $G$ with a generating set $X$. Let
$T^{\alpha}$ be a Schreier transversal for a subgroup $G_1$ of $G$
with respect to $\GA$. Then a {\em Schreier basis} $Y$ for $G_1$ with
respect to $T^{\alpha}$ consists of the non-trivial (in $G$) elements
of the form
\begin{equation}
y_{v,x}=vx(p(vx))^{-1}, \ \ \ \ \ \ \ \ v \in T, \ x \in X,
\end{equation}
where $p$ is the coset map. Notice that the $y_{v,x}$ are not
necessarily distinct elements of $G$.
}
\end{definition}
\begin{proposition}
Let $\GA$ be a presentation of a group $G$ and let $\GAia$ be a presentation
by a Schreier basis of a subgroup $G_1$ of $G$ of finite index. Then
\begin{equation}
\hat \Xi(\GAia) \leq |G:G_1| \hat \Xi(\GA).
\end{equation}
\label{prC87}
\end{proposition}
{\em Proof}.
If $G$ is finite then $\hat \Xi(\GAia) = |G:G_1| \hat \Xi(\GA)$.
Assume that $G$ is infinite. Let $X$ be the set, of cardinality $n$, of
generators of $\GA$, and let $Y$ be the Schreier basis of $\GAia$,
which is of cardinality $\leq 1+ |G:G_1|(n-1)$.
Let ${\mathcal G}, {\mathcal G}_1$ be the Cayley graphs associated with
$\GA,\GAia$
respectively. Let ${\mathcal T}' \in \CF({\mathcal G})$ be the Schreier tree by
which  $Y$ is defined, and let ${\mathcal G}' \in \CF({\mathcal G})$ be the
subgraph induced by ${\mathcal T}'$. Given $\epsilon >0$,
let ${\mathcal G}'_1 \in \CF({\mathcal G}_1)$ satisfy
\begin{equation}
\frac{1-|E^Y_{out}({\mathcal G}'_1)|} {\beta_0({\mathcal G}'_1)} >
\hat \Xi(\GAia) - |G:G_1| \epsilon,
\end{equation}
and let ${\mathcal G}'' \in \FF({\mathcal G})$ be the subgraph induced by
$\bigcup_{g \in V(\scriptstyle {\mathcal G}'_1)} g {\mathcal G}'$.
The edges $E^X_{out}({\mathcal G}'')$ are in 1-1 correspondence with
the edges $E^Y_{out}({\mathcal G}'_1)$. Then
\begin{equation}
\frac{1 - |E^X_{out}({\mathcal G}'')|} {\beta_0({\mathcal G}'')} =
\frac{1- |E^Y_{out}({\mathcal G}'_1)|} {|G:G_1| \beta_0({\mathcal G}'_1)}
> \frac{\hat \Xi(\GAia)}{|G:G_1|} - \epsilon.
\end{equation}
Thus we showed that
\begin{equation}
\hat \Xi(\GAia) \leq |G:G_1| \hat \Xi(\GA).
\end{equation}
\hfill $\Box$ \\

For example, when $G$ is free then we get an equality in
Proposition~\ref{prC87}.
%, but this need not be the case in general.
\section{Direct Products}
\begin{theorem}
Let $\GAi$ be presentations of non-trivial groups $G_i$, for $i=1,2$,
with disjoint generating sets of cardinalities $n_i$ respectively,
and let $\GA$ be the induced presentation of $G=G_1 \times G_2$. Then
\begin{equation}
\hat \Xi(\GA) = \hat \Xi(\GAia) + \hat \Xi(\GAib) - \left( \frac{1}{|G_1|} +
\frac{1}{|G_2|} - \frac{1}{|G|} \right).
\label{eqC88}
\end{equation}
\label{thC20}
\end{theorem}
{\em Proof}.
The claim is true when both $G_1$ and $G_2$ are finite, because then
\begin{equation}
\hat \Xi(\GA) = \frac{1}{|G|}.
\end{equation}

When at least one of the groups is infinite we will show first that
(\ref{eqC88}) is an upper bound for $\hat \Xi(\GA)$.
Suppose that ${\mathcal G}' \in \CFM({\mathcal G})$, where ${\mathcal G}$ is
the Cayley graph of $G$. Then ${\mathcal G}'$ is the union of the
subgraphs ${\mathcal H}'_1, {\mathcal H}'_2$, where ${\mathcal H}'_i$, $i=1,2$,
is the
subgraph generated by the edges with labels in $X_i \cup X_i^{-1}$.
We may also assume that both ${\mathcal H}'_i$ are not empty, because otherwise
we can obtain at most $\max(\hat \Xi(\GAia)-n_2, \hat \Xi(\GAib)-n_1)$
which is
less than or equals the right hand side of (\ref{eqC88}). Then
%For each $i$ let ${\mathcal H}_{i,j}$, $j=1, \ldots, k_i$, be the (non-trivial)
%components of ${\mathcal H}'_i$. Each such component ${\mathcal H}_{i,j}$ may b
e
%regarded as embedded in the Cayley graph ${\mathcal G}_i$ of $G_i$. Then
\begin{equation}
\mu({\mathcal G}') - (n_1+n_2) \leq
(\mu({\mathcal H}'_1) -n_1 ) +
\left( \frac{\beta_1({\mathcal H}'_2)}{\beta_0({\mathcal H}'_2)} -n_2 \right).
%&\leq& \left( \frac{\beta_1({\mathcal H}'_1)+k_1}{\beta_0({\mathcal H}'_1)}
%-n_1 \right) +
%\left( \frac{\beta_1({\mathcal H}'_2)}{\beta_0({\mathcal H}'_2)} -n_2 \right).
\label{eqC92}
\end{equation}
By (\ref{eqC82})
\begin{equation}
\mu({\mathcal H}'_1) -n_1 \leq \hat \Xi(\GAia),
\end{equation}
and
\begin{equation}
\frac{\beta_1({\mathcal H}'_2)}{\beta_0({\mathcal H}'_2)} -n_2
\leq \min(\hat \Xi(\GAib),0).
\end{equation}
Therefore, by symmetry, we get from (\ref{eqC92}) that
\begin{equation}
\hat \Xi(\GA) \leq \min(\hat \Xi(\GAia),\hat \Xi(\GAib),\hat \Xi(\GAia)+
\hat \Xi(\GAib)).
\end{equation}
But these bounds are exactly the ones that appear in (\ref{eqC88}) in case
at least one of the groups is infinite.

It remains to show that (\ref{eqC88}) is really achieved.
Given $\epsilon >0$ then for $i=1,2$ let ${\mathcal H}_i \in \CF({\mathcal G}_i)
$,
where ${\mathcal G}_i$ is the Cayley graph of $\GAi$, satisfying
\begin{equation}
\mu({\mathcal H}_i) -n_i > \hat \Xi(\GAi) - \frac{\epsilon}{4}.
\end{equation}
Suppose also that if $G_i$ is finite then ${\mathcal H}_i = {\mathcal G}_i$
and otherwise $\beta_0({\mathcal H}_i) > 4/\epsilon$.
Let ${\mathcal G}'$ be the subgraph of ${\mathcal G}$ which is the
cartesian product ${\mathcal H}_1 \times {\mathcal H}_2$. Then
\begin{eqnarray}
\lefteqn{\mu({\mathcal G}') - (n_1+n_2) =
 \frac{\beta_1({\mathcal H}_1)\beta_0({\mathcal H}_2) +
\beta_1({\mathcal H}_2)\beta_0({\mathcal H}_1) + 1}
{\beta_0({\mathcal H}_1)\beta_0({\mathcal H}_2)} - (n_1 + n_2)} \\
& & = (\mu({\mathcal H}_1) - n_1 ) + (\mu({\mathcal H}_2) - n_2 )
- \left(\frac{1}{\beta_0({\mathcal H}_1)} + \frac{1}{\beta_0({\mathcal H}_2)} -
\frac{1}{\beta_0({\mathcal H}_1)\beta_0({\mathcal H}_2)} \right) \nonumber \\
& & >  \hat \Xi(\GAia) + \hat \Xi(\GAib) - \left( \frac{1}{|G_1|} +
\frac{1}{|G_2|} - \frac{1}{|G|} \right) - \epsilon,
\nonumber
\end{eqnarray}
and the proof is complete.
\hfill $\Box$ \\
\section{Free Products}
\label{secFP}
Let us look at what happens with the computation of $\hat \Xi(\GA)$ for free
products. We recall that the decomposition of a group $G$ into non-trivial
freely indecomposable factors is unique, up to isomorphism of the factors,
as follows by the Kurosh Subgroup Theorem. Such a decomposition contains
finitely many factors when $G$ is of finite rank, and in fact, by
the corollary to the Grushko-Neumann
Theorem, if $G=G_1*G_2$ then $\mbox{rank}(G)=\mbox{rank}(G_1)+
\mbox{rank}(G_2)$ (see \cite{Lyn}, p. 178).
The following expression plays an important role in computing the value
of $\hat \Xi(\GA)$.
\begin{definition}
{\em
Let $\GA$ be an $n$-generated presentation of a non-trivial group $G$.
If ${\mathcal G}' \in \FFM({\mathcal G})$ let
\begin{equation}
\psi({\mathcal G}') = \frac{\beta_2({\mathcal G}')}{\beta_0({\mathcal G}')-1}.
\end{equation}
Then we define
\begin{equation}
\Psi(\GA) = \sup_{{\mathcal G}' \in \stFFM({\mathcal G})}\psi({\mathcal G}').
\label{eqC90}
\end{equation}
and
\begin{equation}
\hat \Psi(\GA) = 1-n+ \Psi(\GA).
\label{eqC90a}
\end{equation}
}
\end{definition}

We call a presentation {\em reduced} if none
of its generators equals the identity element in the group. Since removing
such ``redundant generators'' does not change the value of $\hat \Xi(\GA)$,
no loss of generality is caused when assuming (as we do in
Theorem~\ref{thC30}), that the presentations are reduced. We say that
a presentation $\GA=<X \mid R>$ is {\em minimal} if
for every proper subset $X'$ of  $X$, $Gp(X') \neq G$. (Here $Gp(X')$ is the
subgroup of $G$ generated by $X'$.)
\begin{theorem}
Let $\GA=<X \mid R>$ be an $n$-generated minimal presentation of
a non-trivial group $G$. Then $1-n \leq \hat \Psi(\GA) \leq 1$.
Moreover,
\begin{description}
\item[(i)] if $G$ is finite then $\hat \Psi(\GA) = n / (|G|-1)$,
with $\hat \Psi(\GA)=1$ if and only if $|G|=2$;
\item[(ii)] if $G$ is infinite then $\hat \Psi(\GA) \leq 0$, with
$\hat \Psi(\GA)= 0$ if and only if $G$ is amenable or $\GA$ is
$2$-generated and one of the generators is of order $2$.
$\hat \Psi(\GA)= 1-n$ if and only if $G$ is free of rank $n \geq 2$. \\
If $H=Gp(Y)$, $Y \subseteq X$, satisfying
\begin{equation}
|H| =
\min \{ |Gp(X')| \mid X' \subseteq X,
|X'|= \max\{ |X''| \mid X'' \subseteq X, |Gp(X'')| < \infty \} \}
\end{equation}
then
\begin{equation}
\Psi(\GA)= \max ( \Xi(\GA), \Psi(\HA)),
\end{equation}
where $\Psi(\HA)$ is calculated as in {\em (i)}.
\end{description}
\label{thC120}
\end{theorem}
{\em Proof}.
Let ${\mathcal G}$ be the Cayley graph of $\GA$.
\begin{equation}
\hat \Psi(\GA) =  1-n+ \sup_{{\mathcal G}' \in \stCFM({\mathcal G})}
\frac{1+(n-1)\beta_0({\mathcal G}') - |E^X_{out}({\mathcal G}')|}
{\beta_0({\mathcal G}')-1} = \sup_{{\mathcal G}'}\frac{n -
|E^X_{out}({\mathcal G}')|} {\beta_0({\mathcal G}')-1}.
\label{eqC103}
\end{equation}
Therefore, $\hat \Psi(\GA) \leq 1$ if and only if
\begin{equation}
|E^X_{out}({\mathcal G}')| \geq n+1-\beta_0({\mathcal G}')
\end{equation}
for every ${\mathcal G}' \in \CFM({\mathcal G})$.
Assume that $X=\{x_1, \ldots, x_n \}$ and $X' = \{x_1, \ldots, x_k \}$,
$0 \leq k \leq n$, is the set of the labels of the edges {\em not}
going-out of ${\mathcal G}'$. If $k=0$ then $|E^X_{out}({\mathcal G}')| \geq n$
and $\hat \Psi(\GA) \leq 0$. Otherwise, ${\mathcal G}'$ is the union of a finite
number of left cosets of $Gp(X')$. Since the presentation is minimal
we have
\begin{equation}
\beta_0({\mathcal G}') \geq |Gp(X')| \geq 2^k.
\end{equation}
Hence
\begin{equation}
|E^X_{out}({\mathcal G}')| \geq n-k \geq n-(2^k-1) \geq n+1-\beta_0({\mathcal G}
').
\label{eqC93}
\end{equation}

Suppose that $G$ is finite. We will show that $\psi$ achieves its maximum on
$\mathcal G$.
Let ${\mathcal G}' \in \CFM({\mathcal G})$ satisfy $\psi({\mathcal G}') \geq
\psi({\mathcal H})$
for every ${\mathcal H} \in \CFM({\mathcal G})$ contained in
${\mathcal G}'$. Let ${\mathcal G}'' \neq {\mathcal G}'$
be a left translate of ${\mathcal G}'$ such that
$V({\mathcal G}' \cap {\mathcal G}'')$ is not empty, and
let ${\mathcal H}' = {\mathcal G}' \cup {\mathcal G}''$. Then
\begin{equation}
\psi({\mathcal H}') \geq \frac{2 \beta_2({\mathcal G}')-\beta_2({\mathcal G}'
\cap {\mathcal G}'')}
{2 \beta_0({\mathcal G}') - \beta_0({\mathcal G}' \cap {\mathcal G}'') -1} \geq
\psi({\mathcal G}'),
\label{eqC94}
\end{equation}
where the right inequality comes from
\begin{equation}
\frac{a}{b-1} \geq \frac{c}{d-1} \ \Llrw \ \frac{2a-c}{2b-d-1} \geq \frac{a}{b-1
},
\end{equation}
with all denominators positive (the case where ${\mathcal G}'$ meets
${\mathcal G}''$ in a single vertex leads to an equality in (\ref{eqC94})).

It is left to examine the case where for every translate ${\mathcal G}''$
of ${\mathcal G}'$, $V({\mathcal G}' \cap {\mathcal G}'')$ is empty.
This means that no edge going-out of
${\mathcal G}'$ has the same label as that of an edge of ${\mathcal G}'$.
Thus ${\mathcal G}'$ is isomorphic to the Cayley graph of $H=Gp(X')$, $|X'|=k$,
and by (\ref{eqC103})
\begin{equation}
\psi({\mathcal G}') = k-1+\frac{k}{|H|-1} \leq k.
\end{equation}
If $H$ is then a proper subgroup of $G$ and $x_j \notin X'$ then the
subgraph ${\mathcal H}'$ which is isomorphic to the Cayley graph of
$H'=Gp(X' \cup X_j)$ satisfies
\begin{equation}
\psi({\mathcal H}') = k+\frac{k+1}{|H'|-1} > k \geq \psi({\mathcal G}').
\end{equation}
We have shown that $\psi({\mathcal G}') \leq \psi({\mathcal G})$ for
every subgraph ${\mathcal G}' \in \CFM({\mathcal G})$. Hence
\begin{equation}
\hat \Psi(\GA) = 1-n+\psi({\mathcal G}) = \frac{n}{|G|-1}.
\end{equation}
Since by the minimality of the presentation $|G| \geq 2^n$ then
$\hat \Psi(\GA)=1$ if and only if $G$ is of order $2$.

Suppose now that $G$ is infinite. If $\psi({\mathcal G}')$ does not have a
maximum on $\CFM({\mathcal G})$ then $\hat \Psi(\GA)=\hat \Xi(\GA)$ since
$|\psi({\mathcal G}')-\xi({\mathcal G}')| \rightarrow 0$
as $\beta_0({\mathcal G}')
\rightarrow \infty$.
The same is true when $\psi({\mathcal G}')$ does have a maximum but there is
no bound to the size of ${\mathcal G}'$ on which $\psi$ achieves its maximum,
e.g. when $G=H*H$ and $H$ is finite. In fact, by Corollary~\ref{crC40}
(vii), when $G=H*H$ then $\Psi(\GA) = \Psi(\HA)$, and if ${\mathcal H}$ is the
Cayley graph of $H$, embedded in ${\mathcal G}$, then when we adjoin $m$
translates of ${\mathcal H}$ to form ${\mathcal K} \in \CFM({\mathcal G})$,
each translate
intersecting the previous one in a single vertex, then
\begin{equation}
\psi({\mathcal K}) = \frac{m \beta_2({\mathcal H})} {m \beta_0({\mathcal H})
-(m-1) -1} = \psi({\mathcal H}).
\end{equation}
When none of the above occurs then $\psi$ achieves its maximum on some
${\mathcal G}' \in \CFM({\mathcal G})$, which is isomorphic to the Cayley graph
of
$H=Gp(X')$, $|X'|=k<n$. Then by (\ref{eqC103})
\begin{equation}
\hat \Psi(\GA)= 1-n+\psi({\mathcal G}') =1-n+( k-1+\frac{k}{|H|-1}) =
k-n+\frac{k}{|H|-1} \leq 0.
\end{equation}
In fact, we see that when $G$ is non-amenable then $\hat \Psi(\GA)= 0$
if and only if $\GA$ is $2$-generated and one of the generators is of
order $2$. When $G$ is infinite amenable then $\hat \Psi(\GA)=0$.

We conclude that for both finite and infinite groups $G$, if
$H=Gp(Y)$, $Y \subseteq X$, satisfies
\begin{equation}
|H| = \min \{ |Gp(X')| \mid X' \subseteq X,
|X'|= \max\{ |X''| \mid X'' \subseteq X, |Gp(X'')| < \infty \} \}
\end{equation}
then
\begin{equation}
\Psi(\GA)= \max ( \Xi(\GA), \Psi(\HA)),
\end{equation}
where $\Psi(\HA)=m-1+m/(|H|-1)$, with $m$ being the number of generators
of $\HA$.

Finally, it is clear that $\hat \Xi(\GA)= 1-n$ if and only if $G$ is
free of rank $n \geq 2$.
\hfill $\Box$ \\

When the presentation is not minimal the assertions of
Theorem~\ref{thC120}
do not hold. For example, if $\GA=< x_1, \ldots, x_n \mid x_1 = x_2 =
\cdots = x_n, \ x_1^2=1>$ and $n \geq 2$ then $\Psi(\GA) = 2n-1 > n$.
\begin{theorem}
For each $i$, $1 \leq i \leq r$, $r \geq 2$, let $\GAi=<X_i \mid R_i>$ be
a reduced $n_i$-generated presentation of a non-trivial group $G_i$
whose Cayley graph is ${\mathcal G}_i$. Let $\GA=<\bigcup_{i=1}^{r} X_i \mid
\bigcup_{i=1}^{r} R_i>$ be the induced $n=\sum_{i=1}^{r} n_i$-generated
presentation of $G=G_1 * G_2 * \cdots * G_r$. Assume also, without loss
of generality, that $\Psi(\GAia) \geq \Psi(\GAib) \geq \cdots \geq
\Psi(\GAir)$, and let $G_{1,2}^{\alpha}$ be the induced presentation of
$G_1*G_2$.
\begin{description}
\item[(i)] If $\Psi(\GAia)=\Xi(\GAia)$ then $\hat \Xi(\GA) = 1-n+\Xi(\GAia)$.
\item[(ii)] If $\Psi(\GAia)=\Psi(\HAia) > \Xi(\GAia)$, where $H_1<G_1$ is
a finite subgroup generated by $Y_1 \subseteq X_1$ as in
Theorem~\ref{thC120}, then
\begin{equation}
\hat \Xi(\GA) =
1-n+ \Xi(G_{1,2}^{\alpha})=
1-n+\max \left( \Xi(\GAia), \Xi(\HAia)+ \frac{\Psi(\GAib)}{|H_1|} \right).
\end{equation}
%\begin{equation}
%\hat \Xi(\GA) = 1-n + \sup_{{\mathcal H} \in \stCF({\mathcal G}_1)}
%\left( \xi({\mathcal H}) + \frac{\Psi(\GAib)}{\beta_0({\mathcal H})} \right)=
%1-n+ \Xi(G_{1,2}^{\alpha}).
%\label{eqC115}
%\end{equation}
\end{description}
\label{thC30}
\end{theorem}
{\em Proof}.
Let ${\mathcal G}$ be the Cayley graph corresponding to $\GA$.
Given $\epsilon > 0$ let ${\mathcal G}' \in \CFM({\mathcal G})$ satisfying
\begin{equation}
\xi({\mathcal G}') > \Xi(\GA) -\epsilon.
\end{equation}
${\mathcal G}'$ has the form
${\mathcal G}' = \bigcup_i {\mathcal H}_i$, where each ${\mathcal H}_i$
is the disjoint union of $k_i \geq 0$ subgraphs ${\mathcal H}_{i,j} \in
\CFM({\mathcal G}_i)$, $j=1, \ldots, k_i$,
which we look at as being subgraphs of ${\mathcal G}_i$, the Cayley graphs
of $\GAi$.
We say that such a subgraph ${\mathcal H}_{i,j}$ is of type $i$.
Starting from some ${\mathcal H}_{i_0,j_0}$, ${\mathcal G}'$ can be constructed
inductively, forming a tree-like structure, by adding at each stage one
of the subgraphs ${\mathcal H}_{i,j}$, which meets the subgraph constructed
up to that stage at a
single vertex, since there are no other simple circuits except the ones
in the subgraphs ${\mathcal H}_{i,j}$ (this is where $\psi$ comes into the
picture). This implies that
\begin{equation}
\xi({\mathcal G}') =
\frac{\sum_{i=1}^{r} \beta_2({\mathcal H}_i)}
{\sum_{i=1}^{r}\beta_0({\mathcal H}_i) +1-\sum_{i=1}^{r}k_i}
= \frac{\sum_{i=1}^{r}\sum_{j=1}^{k_i}\beta_2({\mathcal H}_{i,j})}
{1+\sum_{i=1}^{r}\sum_{j=1}^{k_i}(\beta_0({\mathcal H}_{i,j})-1)}.
\label{eqC95}
\end{equation}
By the form of (\ref{eqC95}) we see that an upper bound for $\Xi(\GA)$ is
$\Psi(\GAia)$, and by using only subgraphs of type 1 and 2 we get a lower
bound $\Xi(\GA) \geq \Psi(\GAib)$. In case $\Psi(\GAia)=\Psi(\GAib)$ then
$\Xi(\GA)=\Psi(\GAia)$. We may then assume that $\GAia$ is involved in
${\mathcal G}'$ when $\epsilon$ is small enough.
We may also assume that except from $\GAia$, ${\mathcal G}'$ involves edges from
other $\GAi$ (otherwise we have a connected subgraph ${\mathcal H}$ of
${\mathcal G}_1$
and we can take two copies of it joined by an edge from some $X_i$,
$i \neq 1$, so that $\xi({\mathcal H})$ is not changed).
We look at the decomposition of ${\mathcal G}'$ into the subgraphs
${\mathcal H}_{i,j}$ as above. Then we reconstruct
$\xi({\mathcal G}')$ in the following way.
We start from a subgraph ${\mathcal H}_{1,j_0}$ of type $1$. Then we add the
subgraph of ${\mathcal G}'$ consisting of some ${\mathcal H}_{i_1,j_1}$,
of type $i_1
\neq 1$ and all the new subgraphs (not including the one we started with)
of type $1$
joined to it (and there may be none of them). We Continue in an
inductive way, so that at the $m$-th stage we add some new
${\mathcal H}_{i_m,j_m}$,
of type $i_m \neq 1$, which is joined to the part constructed up to that
stage, and all the new subgraphs of type $1$ joined to ${\mathcal H}_{i_m,j_m}$.
We finish after
we cover the whole of ${\mathcal G}'$. We show now that there exists
a subgraph ${\mathcal G}''$ of ${\mathcal G}$, decomposed into copies of only tw
o
subgraphs ${\mathcal H} \subseteq {\mathcal G}_1$ and ${\mathcal H}' \subseteq
{\mathcal G}_2$,
such that $\xi({\mathcal G}'') \geq \xi({\mathcal G}')$. First we notice that
if for
some ${\mathcal H}_{i_m,j_m}$, $i_m \neq 1$, the number of subgraphs of type
$1$ joined to it is less than $\beta_0({\mathcal H}_{i_m,j_m})$
then there are subgraphs of type $1$ that we can add to it so that the
$\xi$ does not decrease since $\Psi(\GAia)
\geq \Psi(\GAi)$ for every $i$. So let us suppose that indeed each
such ${\mathcal H}_{i_m,j_m}$ is joined to $\beta_0({\mathcal H}_{i_m,j_m})$
subgraphs of type $1$. Let ${\mathcal H}_{1,j_{k}}$, $k=1, \ldots, p_m=
\beta_0({\mathcal H}_{i_m,j_m})-1$, be the new subgraphs of type $1$
added at the $m$-th stage in the reconstruction of ${\mathcal G}'$.
That is, at that stage $\beta_2$ is increased by
\begin{equation}
a_m= \beta_2({\mathcal H}_{i_m,j_m})+ \sum_{k=1}^{p_m}
\beta_2({\mathcal H}_{1,j_{k}})
\end{equation}
and $\beta_0$ is increased by
\begin{equation}
b_m=\sum_{k=1}^{p_m}\beta_0({\mathcal H}_{1,j_{k}}).
\end{equation}
Since $\Psi(\GAib) \geq \Psi(\GAi)$ for every $i > 2$,
then there exist finite connected subgraphs ${\mathcal H}'_m $ of type
$2$ and ${\mathcal H}''_m$ of type $1$, such that
\begin{eqnarray}
\xi({\mathcal H}''_m) + \frac{\psi({\mathcal H}'_m)}{\beta_0({\mathcal H}''_m)}
&=&
\frac{(\beta_0({\mathcal H}'_m)-1)\beta_2({\mathcal H}''_m) +
\beta_2({\mathcal H}'_m)}
{(\beta_0({\mathcal H}'_m)-1)\beta_0({\mathcal H}''_m)} \\
&\geq&
\frac{(\beta_0({\mathcal H}_{i_m,j_m})-1)\beta_2({\mathcal H}_{1,j_{k}}) +
\beta_2({\mathcal H}_{i_m,j_m})}
{(\beta_0({\mathcal H}_{i_m,j_m})-1)\beta_0({\mathcal H}_{1,j_{k}})} \nonumber
\end{eqnarray}
for every $1 \leq m \leq t=\sum_{i=2}^{r}k_i$ and $1 \leq k \leq p_m$
(first choose ${\mathcal H}'_m$ such that $\psi({\mathcal H}'_m) \geq
\psi(H_{i_m,j_m})$
for every $m$, and then an appropriate ${\mathcal H}''_m$).
Hence it follows that
\begin{equation}
\xi({\mathcal H}''_m) + \frac{\psi({\mathcal H}'_m)}{\beta_0({\mathcal H}''_m)}
\geq \frac{a_m}{b_m}
\end{equation}
(as seen after clearing denominators).
Let ${\mathcal H} \in \CF({\mathcal G}_1), {\mathcal H}' \in \CF({\mathcal G}_2)
$
such that
\begin{equation}
\xi({\mathcal H}) + \frac{\psi({\mathcal H}')}{\beta_0({\mathcal H})}
\geq \xi({\mathcal H}''_m) + \frac{\psi({\mathcal H}'_m)}{\beta_0({\mathcal H}''
_m)}
\end{equation}
for each $m=1, \ldots, t$. By the last two inequalities
the subgraph ${\mathcal G}''$ constructed by
starting with ${\mathcal H}_{1,j_0}$ and adjoining $t$ times the subgraph
consisting of a copy of ${\mathcal H}'$ and $(\beta_0({\mathcal H}') -1)$ copies
of ${\mathcal H}$ satisfies
\begin{equation}
\xi({\mathcal G}'') \geq \xi({\mathcal G}').
\end{equation}
Since $\epsilon$ was chosen arbitrarily, we get by the form of ${\mathcal G}''$
that
\begin{eqnarray}
\label{eqC135}
\hat \Xi(\GA) &=& 1-n + \sup_{{\mathcal H}'_1 \in \stCF({\mathcal G}_1),
{\mathcal H}'_2 \in \stCF({\mathcal G}_2)}
\left( \xi({\mathcal H}'_1) + \frac{\psi({\mathcal H}'_2)}{\beta_0({\mathcal H}'
_1)}
\right) \\
&=& 1-n + \sup_{{\mathcal H}'_1 \in \stCF({\mathcal G}_1)}
\left( \xi({\mathcal H}'_1) + \frac{\Psi(\GAib)}{\beta_0({\mathcal H}'_1)} \right)
\nonumber \\
&=& 1-n+ \Xi(G_{1,2}^{\alpha})
\nonumber
\end{eqnarray}
Let us define $\zeta({\mathcal H}'_1)$ on $\CF({\mathcal G}_1)$
by
\begin{equation}
\zeta({\mathcal H}'_1) = \xi({\mathcal H}'_1) + \frac{\Psi(\GAib)}
{\beta_0({\mathcal H}'_1)}.
\end{equation}
Thus we need to find
\begin{equation}
\Xi(\GA)= \sup_{{\mathcal H}'_1 \in \stCF({\mathcal G}_1)} \zeta({\mathcal H}'_1
).
\end{equation}
If $\Psi(\GAia) = \Psi(\GAib)$ then as we have seen
$\Xi(\GA)=\Psi(\GAia) =\zeta(1)$, where $1$
is the trivial subgraph. Also, when $\Xi(\GAia)=\Psi(\GAia)$ then $\Xi(\GA) =
\Xi(\GAia)$ because $\Psi(\GA) = \Psi(\GAia)$ by (\ref{eqC95}) and
$\Xi(\GA) \geq \Xi(\GAia)$ by the embedding of ${\mathcal G}_1$ in
${\mathcal G}$. When $\Xi(\GAia)<\Psi(\GAia)$ but an arbitrarily large subgraph
${\mathcal H}'_1$ can be chosen without decreasing $\zeta({\mathcal H}'_1)$
then we get that $\Xi(\GA)=\Xi(\GAia)$ as the second summand in the
expression defining $\zeta({\mathcal H}'_1)$ tends to zero when
$\beta_0({\mathcal H}'_1) \rightarrow \infty$.
It remains to check the case where $\zeta$
achieves its maximum on a finite number of members of $\CF({\mathcal G}_1)$.
Let ${\mathcal H}'_1$ be maximal (with respect to the number of vertices)
among these subgraphs. We will show that ${\mathcal H}'_1$ is isomorphic
to the Cayley graph of a finite subgroup $H_1$
of $G_1$ of the form described in
Theorem~\ref{thC120}~(ii) on which $\psi$
achieves its maximum.
For we are given that ${\mathcal H}'_1$ satisfies
$\zeta({\mathcal H}'_1) \geq \zeta({\mathcal H}')$ for every ${\mathcal H}'
\in \CF({\mathcal G})$ which is contained in ${\mathcal H}'_1$.
Suppose there exists ${\mathcal H}''_1 \neq {\mathcal H}'_1$ a left translate of
${\mathcal H}'_1$ such that
${\mathcal H}'_1 \cap {\mathcal H}'' \neq \emptyset$.
Let ${\mathcal H}_1 = {\mathcal H}'_1 \cup {\mathcal H}''_1$. Then
\begin{equation}
\zeta({\mathcal H}_1) \geq \frac{2 \beta_2({\mathcal H}'_1)-\beta_2({\mathcal H}
'_1
\cap {\mathcal H}''_1)+\Psi(\GAib)}
{2 \beta_0({\mathcal H}'_1) - \beta_0({\mathcal H}'_1 \cap {\mathcal H}''_1)} \geq
\zeta({\mathcal H}'_1),
\end{equation}
where the right inequality comes from
\begin{equation}
\frac{a+\Psi}{b} \geq \frac{c+\Psi}{d} \ \Llrw \ \frac{2a-c+\Psi}{2b-d} \geq
\frac{a+\Psi}{b},
\end{equation}
with all denominators positive. But this contradicts the maximality of
${\mathcal H}'_1$. Thus the set of labels of the edges of ${\mathcal H}'_1$
is disjoint from the set of
labels of its outer edges. This means that ${\mathcal H}'_1$ is isomorphic
to the Cayley graph of a subgroup of $G_1$. We need to show that
$\psi$ too achieves its maximum on ${\mathcal H}'_1$. Recall that we are in the
case where $\Psi(\GAia) > \Xi(\GAia)$, and so
$\psi$ achieves its maximum on some subgraph
${\mathcal H}_1$ which is isomorphic to
the Cayley graph of a finite subgroup $H_1 < G_1$. Let ${\mathcal K}_1
\in \CFM({\mathcal G}_1)$ be isomorphic to the Cayley graph of a finite
subgroup $K_1 < G_1$. Thus $\Psi(\HAia) \geq \Psi(K_{1}^{\alpha_1})$.
If $|H_1| \leq |K_1|$ then since for finite groups
\begin{equation}
\Psi(\HAia) \geq \Psi(K_1^{\alpha_1}) \ \Llrw \
\Xi(\HAia) \geq \Xi(K_1^{\alpha_1}),
\end{equation}
we get that
\begin{equation}
\zeta({\mathcal H}_1)=\Xi(\HAia)+\frac{\Psi(\GAib)}{|H_1|} \geq
\Xi(K_1^{\alpha_1})+\frac{\Psi(\GAib)}{|K_1|}=\zeta({\mathcal K}_1).
\end{equation}
So assume $|H_1| > |K_1|$. Then $\psi({\mathcal H}_1) \geq
\psi({\mathcal K}_1)$ is equivalent to
\begin{equation}
\xi({\mathcal H}_1)+\frac{\psi({\mathcal H}_1)}{\beta_0({\mathcal H}_1)} \geq
\xi({\mathcal K}_1)+\frac{\psi({\mathcal H}_1)}{\beta_0({\mathcal K}_1)}.
\end{equation}
Since $\psi({\mathcal H}_1) \geq \Psi(\GAib)$ we have
\begin{equation}
\xi({\mathcal H}_1)-\xi({\mathcal K}_1) \geq \psi({\mathcal H}_1)
\left(\frac{1}{\beta_0({\mathcal K}_1)}
-\frac{1}{\beta_0({\mathcal H}_1)} \right) \geq
\Psi(\GAib) \left(\frac{1}{\beta_0({\mathcal K}_1)} -
\frac{1}{\beta_0({\mathcal H}_1)} \right).
\end{equation}
That is
\begin{equation}
\zeta({\mathcal H}_1)=
\xi({\mathcal H}_1)+\frac{\Psi(\GAib)}{\beta_0({\mathcal H}_1)} \geq
\xi({\mathcal K}_1)+\frac{\Psi(\GAib)}{\beta_0({\mathcal K}_1)}=
\zeta({\mathcal K}_1).
\end{equation}

We conclude that if $H_1$ is a finite subgroup of $G_1$ ($G_1$ may be
finite or infinite) generated by some set $Y_1 \subseteq X_1$ satisfying
\begin{equation}
|H_1| =
\min \{ |Gp(X'_1)| \mid X'_1 \subseteq X_1,
|X'_1|= \max\{ |X''_1| \mid X''_1
\subseteq X_1, |Gp(X''_1)| < \infty \} \}
\end{equation}
then
\begin{equation}
\hat \Xi(\GA) = 1-n+\max \left( \Xi(\GAia),
\Xi(\HAia)+ \frac{\Psi(\GAib)}{|H_1|} \right).
\end{equation}
If, on the other hand, $Gp(X'_1)$ is infinite for every $X'_1
\subseteq X_1$ then $\hat \Xi(\GA) = 1-n+\max \Xi(\GAia)$.
\hfill $\Box$ \\
\begin{example}
{\em
Let $G=G_1*G_2$, where $G_1=C_2 * C_3$ ($C_i$ the cyclic group of
order $i$)
and $G_2=C_4$, with all cyclic factors single-generated.
Then $\Psi(G_1)=\Psi(C_2)=1$, $\Psi(G_2)=1/3$,
$\Xi(G_1)= \Xi(C_2)+\Psi(C_3)/
|C_2|= 1/2+1/(2 \cdot 2) = 3/4$,
and $\Xi(G)=\Xi(G_1)=3/4 > 2/3 = 1/2+1/(3 \cdot 2)=\Xi(C_2)
+\Psi(G_2)/|C_2|$.

On the other hand, suppose that $G_1=C_2 *
C_4$, $G_2=C_3$ and $G=G_1*G_2$. Then
$\Psi(G_1)=1$, $\Psi(G_2)=1/2$, $\Xi(G_1)= 1/2+1/(3 \cdot 2)=2/3$
and $\Xi(G)= 3/4 = 1/2+1/(2 \cdot 2)=
1/\Xi(C_2) +\Psi(G_2)/|C_2| > \Xi(G_1)$.
We see that when $\Psi(G_1) > \Xi(G_1)$ then either of
the possibilities for $\Xi(G)$ that are stated in Theorem~\ref{thC30}
can occur.
}
\end{example}

In the following corollary the results either follow immediately from
Theorem~\ref{thC30} or are already stated within the proof of
Theorem~\ref{thC30} or follow from the arguments there. Therefore
only a partial proof is given.
We further assume that the presentations are reduced.
%and that cyclic and free factors $G_i$ are generated by
%$\mbox{rank}(G_i)$ elements.
\begin{corollary}
Let $G=G_1 * G_2 * \cdots * G_r$ with the assumptions of Theorem~\ref{thC30}.
Then the following claims hold.
\begin{description}
\item[(i)] $\Psi(\GAib) \leq \Xi(\GA) \leq \Psi(\GAia)$.
\item[(ii)] $\Xi(\GA) \geq max_i \{ \Xi(\GAi) \}$. \\
$\Xi(\GA) = \Xi(\GAia)$ if
\begin{description}
\item[(a)] $\Xi(\GAia)=\Psi(\GAia)$; or
\item[(b)] $Gp(X'_1)$ is infinite for every $X'_1 \subseteq X_1$, in
particular if $G_1$ is torsion-free; or
\item[(c)] $\xi({\mathcal H})$ does not have a maximum on $\FF({\mathcal G}_1)$
and $(c(\GAia)-1)\Psi(\GAib) \leq 1$; or
\item[(d)] $\Psi(K_1^{\alpha_1}) \geq \Psi(\GAib)$, where $K_1=Gp(Z_1)$,
$Z_1=X_1-Y_1$ is non-empty and $Y_1 \subseteq X_1$ generates a finite
subgroup $H_1$ for which $\Psi(\GAia)=\Psi(\HAia) > \Xi(\GAia)$.
\end{description}
\item[(iii)] $\hat \Xi(\GA) \leq 1-r+ \sum_{i=1}^{r} \hat \Xi(\GAi)$,
with equality
holding if and only if either $G_2$ is free or $G_1 \simeq G_2 \simeq
C_2$ and $G_3$ (if exists) is free.
\item[(iv)] $\hat \Xi(\GA)= n_1 - n$ if $G_1$ is infinite amenable.
\item[(v)] $\Xi(\GA)= \Xi(\GAia) + \frac{\Psi(G_2^{\alpha_2})}{|G_1|} =
n_1-1+ \frac{\Psi(G_2^{\alpha_2})+1}{|G_1|}$ if $G_1$ is finite. \\
$\Xi(\GA)= n_1-1+ \frac{n_2 |G_2|} {|G_1|(|G_2|-1)}$ if
$G_1$ and $G_2$ are finite.
\item[(vi)] $\hat \Xi(\GA) = 1-r+\frac{1}{|G_1|(1-\frac{1}{|G_2|})}$ if
all the factors are (finite or infinite) cyclic and single-generated.
\item[(vii)] $\Psi(\GA) = \Psi(\GAia)$.
\end{description}
\label{crC40}
\end{corollary}
{\em Proof}.
%(i) The left inequality follows by taking ${\mathcal H}$ to be trivial in
%(\ref{eqC115}). The right inequality is clear from (\ref{eqC95}).
%
%(ii) The inequality follows from the embedding of each ${\mathcal G}_i$ in
%${\mathcal G}$, so that ${\mathcal G}'$ can be taken as a subgraph of
%${\mathcal G}_i$.
%If $\Xi(\GAia)=\Psi(\GAia)$ (e.g. when $G_1$ is free or infinite amenable)
%then by (i) we get that $\Xi(\GA) = \Xi(\GAia)$.
{\bf (ii) (c)} For every $\epsilon > 0$ there exists
${\mathcal H} \in \CF({\mathcal G}_1)$ such that
\begin{equation}
\Xi(\GA) \leq \zeta({\mathcal H}) + \epsilon =
\xi({\mathcal H}) + \frac{\Psi(\GAib)}{\beta_0({\mathcal H})} +\epsilon.
\end{equation}
When $\xi({\mathcal H})$ does not have a maximum on $\FF({\mathcal G}_1)$ (which
is the case ``in general'') then by Theorem~\ref{thC10}
\begin{equation}
\Xi(\GAia) \geq \xi({\mathcal H})+\frac{1}{(c(\GAia)-1) \beta_0({\mathcal H})}.
\end{equation}
If, in addition, $(c(\GAia)-1)\Psi(\GAib) \leq 1$ then we get from the
two inequalities that
\begin{equation}
\Xi(\GA) \leq \Xi(\GAia)),
\end{equation}
and by the inequality in the other direction this is an equality.

{\bf (ii) (d)} $Gp(Y_1 \cup Z_1) < G_1$ is a quotient of $L_1 =H_1*
K_1$ and therefore, by Theorem~\ref{thC40}
$\Xi(\GAia) \geq \Xi(L_1^{\alpha_1}) \geq
\Xi(\HAia)+ \frac{\Psi(L_1^{\alpha_1})}{|H_1|} \geq
\Xi(\HAia)+ \frac{\Psi(\GAib)}{|H_1|}$. Thus, by Theorem~\ref{thC30},
$\Xi(\GA) = \Xi(\GAia)$.

{\bf (iii)} $\hat \Xi(\GA) \leq 1-r+ \sum_{i=1}^{r} \hat \Xi(\GAi)$ since by
(\ref{eqC95})
\begin{equation}
\xi({\mathcal G}') = \sum_{i=1}^{r} \frac{\beta_2({\mathcal H}_i)}
{\beta_0({\mathcal G}')}
\leq \sum_{i=1}^{r} \xi({\mathcal H}_i).
\end{equation}
A necessary condition for equality in (iii) is that for every $\epsilon >0$
there exists ${\mathcal G}'$ such that for {\em every} non-free factor $G_i$,
$\beta_0({\mathcal H}_i) / \beta_0({\mathcal G}') > 1-\epsilon$. But this
is possible if and only if there is either only one non-free factor,
or there are
two such factors and each is isomorphic to $C_2$, the group
of order $2$.

{\bf (iv)} If $G_1$ is infinite amenable then $\Xi(\GAia)=\Psi(\GAia)$
and therefore $\hat \Xi(\GA)= 1-n+\Xi(\GAia)=n_1 - n$.

{\bf (v)} When $G_1$ is finite then $\zeta({\mathcal H}'_1)$, ${\mathcal H}'_1 \
in
\CF({\mathcal G}_1)$ achieves its maximum on ${\mathcal G}_1$.

{\bf (vi)} When $G_i$ is cyclic then $\Xi(\GAi)= (|G_i|)^{-1}$ and $\Psi(\GAi)=
(|G_i|-1)^{-1}$. Then when $G_{1,2}^{\alpha}$ is the
induced presentation of $G_1 * G_2$ we get from Theorem~\ref{thC30} that
\begin{equation}
\hat \Xi(\GA) = 1-r+ \Xi(G_{1,2}^{\alpha}) = 1-r+\frac{1}{|G_1|} +
\frac{1}{(|G_2|-1)|G_1|}=
-1+\frac{1}{|G_1|(1-\frac{1}{|G_2|})}.
\end{equation}

{\bf (vii)} This is also clear by (\ref{eqC95}).
\hfill $\Box$
\section{The Normalized Balanced Cyclomatic Quotient}
The definition we used for $\hat{\Xi}(\GA)$ says that its value has to be
looked for in the ``best chain'' we can find in the graph. If we
look on the other hand only on the chain which consists of the
concentric balls around $1$, we get some ``averaging'' and in this
case the value we calculate depends on the growth of the group.
\begin{definition}
{\em
Let $\GA$ be $n$-generated and let ${\mathcal G}$
be the corresponding Cayley graph. Let ${\mathcal B}_i$ be the (induced)
concentric balls around $1$ of radius $i$ in ${\mathcal G}$.
Then we define the {\em normalized balanced cyclomatic quotient} of $\GA$ by
\begin{equation}
\hat \Theta(\GA) = 1-n + \limsup_{i \rightarrow \infty} \xi({\mathcal B}_i).
\end{equation}
}
\end{definition}

Clearly $\hat \Theta(\GA) \leq \hat \Xi(\GA)$, and equality holds when $G$ has
subexponential growth, as seen from the proposition below.
When $G$ is infinite then $\hat \Theta(\GA)$ takes values between $1-n$
and $0$.
\begin{proposition}
The growth of $G$ is exponential if and only if $\hat \Theta(\GA)<0$.
\label{prC30}
\end{proposition}
{\em Proof}.
This follows by having $1-n+\xi({\mathcal B}_i) =
(1-|E^X_{out}({\mathcal B}_i)|)
/ \beta_0({\mathcal B}_i)$, and $\beta_0({\mathcal B}_i)$, $i=0,1,2,
\ldots$ is (a representative of the equivalent class of) the growth
function of $G$.
\hfill $\Box$ \\

We say that a non-empty subgraph $\mathcal H \subseteq \mathcal G$ has
{\em thickness} $\geq r$ if $\mathcal H$ contains a non-empty subgraph
$\mathcal H'$
such that $d(v,{\mathcal H}')=r$ for every $v \in \partial \mathcal H$.
The supremum
on all such $r$ is the thickness of $\mathcal H$. When $\partial
\mathcal H$ is
empty then the thickness of $\mathcal H$ is $\infty$. Thus every subgraph
has thickness $\geq 0$, and it has thickness $\geq 1$ if and only if
every vertex on its boundary is adjacent to an interior vertex.

We call a subgraph $H < F$ a {\em supnormal} subgroup if the maximal normal
subgroup $N$ of $F$ which is contained in $H$ is non-trivial.
In the next lemma we refer to the length of a shortest word in $N$, which
is the girth of the Cayley graph of $F/H$ in case $H$ is normal.
\begin{lemma}
Let $H$ be a supnormal subgroup of the free group $F$ of rank $n$
and let $m>0$ be the length of a shortest word in the maximal normal
subgroup of $F$ contained in $H$. Then for every finite subgraph
${\mathcal G}'$ of the cosets graph of $H$ with thickness $> m/2$
\begin{equation}
\xi({\mathcal G}') \geq \frac{2(n-1)}
{m((2n-1)^{1+m/2} -1)}.
\end{equation}
\label{lmC10}
\end{lemma}
{\em Proof}.
Clearly it is enough to show the claim holds for connected subgraphs.
Let $\lambda$ be a fixed circuit of length $m$ corresponding to an
element as in the lemma. Let ${\mathcal G}'$ be a connected finite subgraph
of thickness $\geq m/2$ in the cosets graph of $H$, and let ${\mathcal G}''$ be
the subgraph of ${\mathcal G}'$ induced in ${\mathcal G}'$ by the set of vertice
s
of distance $\geq m/2$ from $\partial {\mathcal G}'$. We cover ${\mathcal G}''$
with copies of $\lambda$, so that each new circuit begins in a vertex not
covered yet. Part of the circuits may extend beyond ${\mathcal G}''$, but not
beyond ${\mathcal G}'$. Then
\begin{equation}
\beta_2({\mathcal G}') \geq \frac{1}{m} \beta_0({\mathcal G}'')
\end{equation}
since each circuit was counted at most $m$ times. As for the cardinality
of ${\mathcal G}'$, we have
\begin{equation}
\beta_0({\mathcal G}') \leq \beta_0({\mathcal G}'')+a
\beta_0(\partial {\mathcal G}'') \leq (1+a)\beta_0({\mathcal G}''),
\end{equation}
where
\begin{equation}
a=\sum_{j=1}^{m/2} (2n-1)^j.
\end{equation}
Therefore
\begin{equation}
\xi({\mathcal G}') \geq
\frac{\beta_0({\mathcal G}'')}{m(1+a)\beta_0({\mathcal G}'')} =
\frac{2(n-1)} {m((2n-1)^{1+m/2} -1)}.
\end{equation}
\hfill $\Box$ \\

As an immediate corollary we have
\begin{proposition}
Let $\GA=<X \mid R>$ be a presentation of a non-free group $G$ with
$|X|=n$. Let $m$ be the length of a shortest relator.
Then
\begin{equation}
\hat \Theta(\GA) \geq 1-n+\frac{2(n-1)} {m((2n-1)^{1+m/2} -1)}.
\end{equation}
\label{prC100}
\end{proposition}
\hfill $\Box$ \\
\baselineskip=6pt  %line can be omitted if changed to 12pt##
%\small
%\input{bib}

\end{document}